%Paper: alg-geom/9506009
%From: voloch@math.utexas.edu
%Date: Tue, 6 Jun 95 10:51:04 -0500

%plain TeX
\normalbaselineskip=1.6\normalbaselineskip\normalbaselines
\magnification=1200
%\NoBlackBox
\overfullrule=0pt

\def\pmb#1{\setbox0=\hbox{#1}%
 \kern-.025em\copy0\kern-\wd0
 \kern.05em\copy0\kern-\wd0
 \kern-.025em\raise.0433em\box0 }

\def\F{{\bf F}}
\def\P{{\bf P}}

\def\mod{\hbox{mod}}

\def \bs{\bigskip}

\def \d{\delta}

\def \a{\alpha}

\def \con{\equiv }
\def \e{\epsilon}

\centerline{\bf Curves that change genus can have arbitrarily many rational
points}
\medskip
\centerline{\bf Jos\'e Felipe Voloch}
\bs

Let $K$ be a global field of positive
characteristic $p$. In other words, $K$ is a function field in one
variable over a finite field of characteristic $p$.
Let $C$ be a projective algebraic curve defined over $K$.
One defines the absolute genus of $C$, in the usual way, by
extending the field to the algebraic closure. We also define
the genus of $C$ relative to $K$ to be the integer $g_K$ that makes the
Riemann-Roch formula hold, that is, for any $K$-divisor $D$ of $C$,
of sufficiently large degree, the
dimension, $l(D)$, of the $K$-vector space of functions of $K(C)$ whose
polar divisor is bounded by $D$, is $\deg D + 1 - g_K$. Since $K$ is
not perfect, the relative genus may change under inseparable extensions.
(See e.g [A] or [T]).
We have shown ([V]) that if the genus of $C$ relative to $K$ is different from
the absolute genus of $C$ then $C(K)$ is finite. The proof in [V] can be
easily adapted to give an upper bound for $\sharp C(K)$, which
depends on $C$. The purpose of this note is to give examples of curves $C/K$
with fixed $g_K$ for which $\sharp C(K)$ is arbitrarily large. The motivation
for considering this problem comes from the work of Caporaso et al. [CHM],
where they show that a conjecture of Lang implies that, for a number field
$K$, $\sharp C(K)$ can be
bounded in terms of $g$ and $K$ only for all curves $C/K$ of genus $g \ge 2$.
It is tempting to see this as a way of disproving Lang's conjecture by
constructing curves with many points. It is also natural to consider the
problem over function fields, but $C(K)$ can be infinite if $C$ is isotrivial.
Otherwise $C(K)$ is finite, if $g \ge 2$, by a classical result of Samuel ([S])
and
one can ask if $\sharp C(K)$ can be uniformly bounded. We don't know the answer
to this question for smooth curves (but see [BV] for a strong bound on
$\sharp C(K)$).

\proclaim Theorem. Let $p > 2$ be a prime and $q=p^n$. Consider the curve
$C_n/\F_p(t)$ defined by $x-(t+t^{q+2}+t^{2q+3}+...+t^{(p-2)q+p-1})x^p=y^p$.
$C_n$ has absolute genus zero but has genus relative to $\F_p(t)$ equal to
$(p-1)(p-2)/2$. Furthermore $\sharp C_n(\F_p(t)) \ge p^{2^n/2n}$ and
$\sharp C_n(\F_{p^{2n}}(t)) \ge p^{2^n}$.

{\it Proof:} It is clear that $C_n$ has absolute genus zero. The statement
about the relative genus holds for any curve $x-ax^p=y^p, a \in K \setminus
K^p$
and can be checked by appealing to Tate's analogue of the Hurwitz formula
([T], theorem 2) to the extension of fields $K(x,y,a^{1/p})/K(x,y)$.

We will construct points on $C_n$ whose $x$-coordinate is
of the form $a(t)/(t^{q+1}-1)$, where $a(t) = \sum_{i=0}^{q-1} \a_it^i$.
We will get a point in $C_n/\bar \F_p(t)$ if $(t^{q+1}-1)^{p-1}a(t)-
(t+t^{q+2}+t^{2q+3}+...+t^{(p-2)q+p-1})a(t)^p$ is a $p$-th power. Using that
$(t^{q+1}-1)^{p-1}=\sum_{i=0}^{p-1}t^{(q+1)i}$ and comparing coefficients,
this condition is equivalent to $\a_i = \a_{(i+q)/p}^p, i \con 0 (\mod p),
\a_i = \a_{(i-1)/p}^p, i \con 1 (\mod p), \a_i = 0,$ otherwise.

Consider
the map $\phi(i)$ defined for positive integers $i, i \con 0,1 (\mod p)$
by $\phi(i) = (i+q)/p, i \con 0 (\mod p),(i-1)/p, i \con 1 (\mod p)$.
It has the following alternate description, for $i < q$.
If $i = \sum_{j=0}^{n-1} \e_jp^j, 0 \le \e_j \le p-1,
\phi(i) = \sum_{j=1}^{n-1} \e_jp^{j-1} + \d p^{n-1}$, where $\d=1$ if
$\e_0=0$ and $\d=0$ if $\e_0=1$. It follows that if $\e_j \ne 0,1$ for some
$j$ then $\phi^r(i) \not\equiv 0,1 (\mod p)$ for some $r > 0$. On the other
hand, if $\e_j = 0,1$ for all $j$ then $\phi^r(i)$ is defined for all $r > 0$.
Moreover it is easy to check that, in this case, $\phi^{2n}(i) = i$.

Returning to our $\a_i$'s, we see that $\a_i = 0$ if $\e_j \ne 0,1$ for
some $j$ and that $\a_{\phi(i)}^p=\a_i$ and $\a_i^{p^{2n}}=\a_i$ if
$\e_j = 0,1$ for all $j$. If $\a_i \in \F_p$ this simply means
$\a_{\phi(i)}=\a_i$.
The set of polynomials $a(t) \in \F_p[t]$ satisfying our
conditions form a $\F_p$-vector space and each orbit of $\phi$ contributes
one dimension to it. Since each orbit has at most $2n$ elements and there
are $2^n$ distinct $i = \sum_{j=0}^{n-1} \e_jp^j, \e_j = 0,1$, we obtain
at least $2^n/2n$ orbits, hence the count for $\F_p(t)$. In the case of
$\F_{p^{2n}}(t)$, an orbit of length $r$ contributes an $r$-dimensional
$\F_p$-vector space. Since $r | 2n$, the theorem follows.

{Remarks:} 1. It can be shown, using the methods of [V], that indeed the
points produced in the proof of the theorem are all the rational points
of $C_n$.

2. In the case $p=3$, $C_n$ is a quasi-elliptic fibration over $\P^1$
in the sense of the classification of surfaces ([BM],[L]) and our result
shows that the 3-rank of the group of sections (the "Mordell-Weil" group)
can be arbitrarily large. (See [I1,2])

3. The curves $C_n$ are the members of the family of curves $x-tf(u)x^p=y^p$,
where $f(u) = \sum_{i=0}^{p-2}u^i$, for $u = t^{p^n+1}$. It follows from
the results of [V] that $tf(u)$ is a $p$-th power in $\F_p(t)$ for only
finitely many $u \in \F_p(t)$, so the curve corresponding to a given $u \in
\F_p(t)$ has finitely many points for all but finitely many $u$'s, again by
the results of [V]. Following [CHM] we consider the total space of the family,
that is, the surface $S$ over $\F_p(t)$ defined by $x-tf(u)x^p=y^p$ and, as is
shown in [CHM], the set of rational points of $S$ is Zariski dense, for
otherwise, the theorem above would be violated. Since $S$ is unirational,
it is not surprising that this holds for some extension of $\F_p(t)$, but
since $S$ cannot be covered by $\P^2$ over $\F_p(t)$, it is surprising that
this occurs over $\F_p(t)$. Also, $S$ is of general type for $p \ge 7$, thus
showing that Lang's conjecture on varieties of general type (see [CHM])
cannot be easily transposed to positive characteristic.

\medskip
{\bf Acknowledgements:} The author would like to thank J. Vaaler
for suggesting an improvement in the treatment of the function $\phi$ above
and the NSF (grant DMS-9301157) and the Alfred P. Sloan Foundation
for financial support.

\medskip

\centerline{\bf References.}
\bigskip
\noindent
[A] E. Artin, {\it Algebraic numbers and algebraic functions}, Gordon and
Bleach,New York, 1967.
\medskip
\noindent
[BM] E. Bombieri, D. Mumford, {\it Enriques' classification of surfaces in
characteristic $p$}, I in {\it Global analysis}, D.C. Spencer, S. Iyanaga,
eds., Princeton Univ. Press, 1969, II in {\it Complex Analysis and Algebraic
Geometry}, Cambridge Univ. Press, 1977, III, Inventiones Math. {\bf 35} (1976)
\medskip
\noindent
[BV] A. Buium, J. F. Voloch, {\it Lang's conjecture in characteristic $p$:
an explicit bound}, Compositio Math., to appear.
\medskip
\noindent
[CHM] L. Caporaso, J. Harris, B. Mazur: {\it Uniformity of
rational points,}  preprint, 1994.
\medskip
\noindent
[I1] H. Ito, {\it The Mordell-Weil groups of unirational quasi-elliptic
surfaces
in characteristic 3}, Math. Z., {\bf 211} (1992) 1-40.
\medskip
\noindent
[I2] H. Ito, {\it The Mordell-Weil groups of unirational quasi-elliptic
surfaces
in characteristic 2}, T\^ohoku Math. J., {\bf 46}(1994),221-251.
\medskip
\noindent
[L] W. E. Lang, {\it Quasi-elliptic surfaces in characteristic 3}
Ann. Sci. \'Ec.
Norm. Sup\'er., IV, Ser. {\bf 12},(1979) 473-500.
\medskip
\noindent
[S] P. Samuel, {\it Compl\'ements \`a un article de Hans
Grauert sur la conjecture de Mordell.} I.H.E.S. Publ. Math. no. 29 (1966)
55-62.
\medskip
\noindent
[T] J. T. Tate, {\it Genus change in inseparable extensions of function
fields}, Proc. Amer. Math. Soc. {\bf 3} (1952) 400-406.
\medskip
\noindent
[V] J. F. Voloch, {\it A diophantine problem on algebraic curves over function
fields of positive characteristic}, Bull. Soc. Math. France {\bf 119} (1991)
121-126.
\medskip
\noindent
Dept. of Mathematics, Univ. of Texas, Austin, TX 78712, USA
\smallskip
\noindent
e-mail: voloch@math.utexas.edu

\end